\documentclass{commat}

%%% AUTHOR'S PACKAGES %%%
\usepackage{booktabs}
\usepackage[nameinlink, noabbrev, capitalize]{cleveref}
\usepackage[shortlabels]{enumitem}
\usepackage[cal=esstix, scr=rsfs]{mathalfa}
\usepackage{tikz-cd}
\usetikzlibrary{decorations.pathmorphing,shapes,shapes.geometric,calc}

%%% AUTHOR'S DEFINITIONS %%%
% for defining ellipses by their focal points
\newcommand\ellipsebyfoci[4]{% options, focus pt1, focus pt2, sum
  \path let \p1=(#2), \p2=(#3), \p3=($(\p1)!.5!(\p2)$)
  in \pgfextra{
    \pgfmathsetmacro{\angle}{atan2(\y2-\y1,\x2-\x1)}
    \pgfmathsetmacro{\focal}{veclen(\x2-\x1,\y2-\y1)/2/1cm}
    \pgfmathsetmacro{\lentotcm}{#4/1cm}
    \pgfmathsetmacro{\axeone}{(\lentotcm - 2 * \focal)/2+\focal}
    \pgfmathsetmacro{\axetwo}{sqrt((\lentotcm/2)*(\lentotcm/2)-\focal*\focal}
  }
  (\p3) node[#1,inner sep=0,rotate=\angle,ellipse,minimum width=2*\axeone cm,minimum height=2*\axetwo cm]{};
}

\renewcommand{\leq}{\leqslant}
\renewcommand{\geq}{\geqslant}
\renewcommand{\cal}[1]{{\mathcal{#1}}}
\newcommand{\scr}[1]{{\mathscr{#1}}}
\newcommand{\longhookrightarrow}{\lhook\joinrel\longrightarrow}
\newcommand{\anotherbullet}{\star}
\newcommand{\cover}{\mathscr{U}}
\newcommand{\anothercover}{\mathscr{V}}
\newcommand{\op}{{\mathrm{op}}}
\newcommand{\sheaf}[1]{{#1}^\mathrm{a}}
\newcommand{\constant}[1]{\underline{#1}}
\newcommand{\category}[1]{\mathsf{#1}}
\newcommand{\simplicial}[1]{\mathsf{s}\,{#1}}
\newcommand{\cosimplicial}[1]{\mathsf{c}\,{#1}}
\newcommand{\Space}{{\category{Space}}}
\newcommand{\Set}{{\category{Set}}}
\newcommand{\Sh}{{\category{Sh}}}
\newcommand{\PSh}{{\category{PSh}}}
\newcommand{\sSh}{{\simplicial{\Sh}}}
\newcommand{\sSet}{{\simplicial{\Set}}}
\newcommand{\sSpace}{{\simplicial{\Space}}}
\newcommand{\cSh}{{\cosimplicial{\Sh}}}
\newcommand{\cSet}{{\cosimplicial{\Set}}}
\newcommand{\scartSh}{{\Sh^\mathrm{s.\,cart}}}
\newcommand{\cartSh}{{\Sh^\mathrm{cart}}}
\newcommand{\Cat}{\category{Cat}}
\newcommand{\Lens}[1]{\category{Lens}\,({#1})}
\newcommand{\ModL}{\category{Mod}^\mathrm{L}}
\newcommand{\id}{\mathrm{id}}
\newcommand{\Ker}{\operatorname{Ker}}
\renewcommand{\Im}{\operatorname{Im}}
\DeclareMathOperator{\Hom}{Hom}
\DeclareMathOperator{\End}{End}
\DeclareMathOperator{\const}{const}
\DeclareMathOperator{\eq}{eq}
\DeclareMathOperator{\Spec}{Spec}
\DeclareMathOperator{\Coh}{\mathsf{Coh}}
\DeclareMathOperator{\CCoh}{\mathsf{CCoh}}
\DeclareMathOperator*{\colim}{colim}
\DeclareMathOperator*{\laxholim}{laxholim}
\DeclareMathOperator*{\holim}{holim}

\title{%
    Various notions of (co)simplicial (pre)sheaves
    }

\author{%
    Timothy Hosgood
    }

\affiliation{\address{Stockholm University}}

\abstract{%
    The phrase ``(co)simplicial (pre)sheaf'' can be reasonably interpreted in multiple ways. In this survey we study how some of these interpretations relate to one another. We end by giving some example applications of the most general of these notions.
    }



\VOLUME{32}
\NUMBER{1}
\YEAR{2024}
\firstpage{73}
\DOI{https://doi.org/10.46298/cm.10359}

\begin{paper}

\section{Introduction}

The phrase ``simplicial sheaf'' turns up in many different places in the literature, and it is not immediately clear that it always has the same meaning.
As just one example\footnote{Appealing to the equivalence between vector bundles and locally free sheaves.}, both \cite{TT1986} and \cite{HT2021} use the phrase ``simplicial vector bundle'', but the former means something like ``a cosimplicial vector bundle over a simplicial space'' whereas the latter means something more along the lines of ``a strongly cartesian simplicial object in (the Grothendieck construction of) the category of vector bundles'' (to use terminology which will be introduced in this present paper).

This survey intends to define the various notions of ``(co)simplicial (pre)sheaf'' known to the author, and describe how these notions relate to one another.
Some of this content was already present in the appendices of \cite{Hos2020}, but what we present here is a more complete version.
We finish by giving some example applications of ``simplicial sheaves'' in their various guises, with the hope that the most general notion (that of a \emph{sheaf on a simplicial space}) can prove to be useful outside of the relatively few applications (known to the author) for which it is currently used.
In the interest of readability, there are many things that will not be discussed (generalisations to Grothendieck topologies, hypercovers, descent, and so on).

In \cref{section:sheaves,section:sections-of-sheaves} we give some basic definitions, and compare the resulting objects;
in \cref{section:cosimplicial-diagrams-of-categories} we describe how these objects arise from a relatively involved construction that might be more familiar to algebraic geometers and algebraic topologists;
finally, in \cref{section:applications}, we describe one ``completely formal'' application of these objects in complex geometry, and one ``rather hand-wavy'' application to certain string diagrams in (symmetric) monoidal categories.

We are very thankful to Damien Calaque for his patient ear and wise words, and to Vincent Wang-Maścianica, without whom \cref{subsection:certain-string-diagrams} would not exist.
We would also like to thank the reviewer for their constructive comments.

\subsection{Conventions and notation}

\begin{itemize}
  \item Unless otherwise stated, sheaves and presheaves are $\Set$-valued (but, in general, most of the constructions and results here will hold for sheaves of abelian objects).
  \item We write $\Space$ to mean some fixed category of suitably nice spaces (e.g. finite CW-complexes, compactly generated topological spaces, Noetherian locally ringed spaces, etc.).
  \item Throughout, let $X$ be a space, and $\cover$ an open cover of $X$.
    We denote by $X_\bullet^\cover$ the corresponding Čech nerve.
    We write $X^\op$ to mean the opposite of the poset category of open sets of $X$ (so that a ``topological'' presheaf on $X$ is exactly a functor $X^\op\to\Set$).
  \item We denote by $\Delta$ the simplex category, with \emph{coface} maps $f_p^i\colon[p-1]\to[p]$ and \emph{codegeneracy} maps $s_i^p\colon[p+1]\to[p]$.
  \item Given a category $\cal{C}$, we write $\simplicial{\cal{C}}$ to denote the category $[\Delta^\op,\cal{C}]$ of simplicial objects in $\cal{C}$, and $\cosimplicial{\cal{C}}$ to denote the category $[\Delta,\cal{C}]$ of cosimplicial objects in $\cal{C}$.
  \item We denote sheafification by $\sheaf{(-)}$. Given a category $[\cal{C}^\op,\cal{D}]$ of $\cal{D}$-valued presheaves, we write $\sheaf{[\cal{C}^\op,\cal{D}]}$ to mean the corresponding (sub)category of sheaves.
\end{itemize}

\section{Sheaves}
\label{section:sheaves}

\subsection{Simplicial sheaves}

\begin{definition}
\label{definition:simplicial-sheaf}
  Let $X$ be a space.
  A \emph{simplicial sheaf on $X$} is either of the two following equivalent things:
  \begin{enumerate}[(i)]
    \item\label{definition:simplicial-sheaf-1}
      a simplicial object in the category of $\Set$-valued sheaves on $X$, i.e.\
      \[
        \cal{F}_\bullet \in [\Delta^\op,\sheaf{[X^\op,\Set]}]
      \]
    \item\label{definition:simplicial-sheaf-2}
      a $\sSet$-valued presheaf on $X$ that is a sheaf in each simplicial degree, i.e.\
      \[
        \cal{F}_\bullet \in [X^\op,[\Delta^\op,\Set]]
        \quad\mbox{such that}\quad
        \cal{F}_p \in \sheaf{[X^\op,\Set]},
      \]
  \end{enumerate}
  where the equivalence is given by $[\cal{C},[\cal{D},\cal{E}]]\simeq[\cal{C}\times\cal{D},\cal{E}]\simeq[\cal{D},[\cal{C},\cal{E}]]$.

  We can similarly define a \emph{cosimplicial sheaf on $X$} as a \emph{cosimplicial} object in the category of $\Set$-valued sheaves on $X$, i.e. $\cal{F}^\bullet\in[\Delta,\sheaf{[X^\op,\Set]}]$.

  We denote the resulting category of simplicial sheaves on $X$ by $\sSh(X)$, and the category of cosimplicial sheaves on $X$ by $\cSh(X)$.
\end{definition}

\subsection{Sheaves on a simplicial space}

Throughout, let $X_\bullet\in\simplicial{\Space}$ be a simplicial space.
We start by giving a ``naive'' definition, as found in \cite[Definition~0.C.4]{Gre1980}, before giving the ``succinct'' categorical definition, as suggested to us by D. I. Spivak.

\begin{definition}
\label{definition:naive-definition-of-sheaf-on-simplicial-space}
  A \emph{sheaf on $X_\bullet$} is the data of a sheaf $\cal{F}^p$ on $X_p$ for all $p\in\mathbb{N}$, along with functorial\footnote{Meaning $\cal{F}^\bullet(\id)=\id$, and $\cal{F}^\bullet(\beta\circ\alpha)=\cal{F}^\bullet(\beta)\circ\cal{F}^\bullet(\alpha)$ for all $\alpha\colon[p]\to[q]$ and $\beta\colon[q]\to[r]$ in $\Delta$.} morphisms
  \[
    \cal{F}^\bullet\alpha\colon(X_\bullet\alpha)^*\cal{F}^p \longrightarrow \cal{F}^q
  \]
  for all $\alpha\colon[p]\to[q]$ in $\Delta$.
  A \emph{morphism} $\varphi^\bullet\colon\cal{F}^\bullet\to\cal{G}^\bullet$ between two such objects is the data of a sheaf morphism $\varphi^p\colon\cal{F}^p\to\cal{G}^p$ for all $p\in\mathbb{N}$ such that the diagram
  \[
    \begin{tikzcd}[sep=huge]
      (X_\bullet\alpha)^*\cal{F}^p
        \rar["(X_\bullet\alpha)^*\varphi^p"]
        \dar["\cal{F}^\bullet\alpha",swap]
      & (X_\bullet\alpha)^*\cal{G}^p
        \dar["\cal{G}^\bullet\alpha"]
    \\\cal{F}^q
        \rar["\varphi^q",swap]
      & \cal{G}^q
    \end{tikzcd}
  \]
  commutes for all $\alpha\colon[p]\to[q]$ in $\Delta$.

  We denote the category of sheaves on the simplicial space $X_\bullet$ by $\Sh(X_\bullet)$.
\end{definition}

\begin{remark}
\label{remark:sheaves-on-simplicial-spaces-look-cosimplicial}
  Note that, although these sheaves live on \emph{simplicial} spaces, they themselves look more like \emph{cosimplicial} objects, since they are \emph{covariant} with respect to the simplex category.
  As we will soon justify, however, they really are simplicial objects: the apparent covariance comes from the fact that the space itself is contravariant with respect to the simplex category, and the sheaf on the space is contravariant with respect to the space, and ``\emph{two $\op$s make an $\id$}''.
\end{remark}

\begin{remark}
  We are being rather agnostic about the category in which our sheaves take values, but one important thing to note is that \cref{definition:naive-definition-of-sheaf-on-simplicial-space} involves \emph{pullbacks} of sheaves.
  If we are simply considering sheaves of sets, then this pullback is the sheaf-theoretic pullback, i.e. the inverse image $f^*=f^{-1}$.
  If, however, we wish to consider ``algebraic'' sheaves (such as the case where $X_\bullet$ is a simplicial \emph{locally ringed} space, and each $\cal{F}^p$ is a sheaf \emph{of $\cal{O}_{X_p}$-modules} on $X_p$) then we need to use the ``algebraic'' pullback
  \[
    f^*(-) = f^{-1}(-)\otimes_{f^{-1}\cal{O}_Y}\cal{O}_X
  \]
  in \cref{definition:naive-definition-of-sheaf-on-simplicial-space}.
\end{remark}

\begin{definition}
  We say that a sheaf $\cal{F}^\bullet$ on $X_\bullet$ is \emph{strictly cartesian}\footnote{\cite{Hos2020a} uses the terminology ``\emph{strongly} cartesian''.} if the $\cal{F}^\bullet\alpha$ are isomorphisms for all $\alpha\colon[p]\to[q]$ in $\Delta$, and denote the (sub)category of such objects by $\scartSh(X_\bullet)$;
  we say that it is (\emph{weakly}) \emph{cartesian} if the $\cal{F}^\bullet\alpha$ are all weak equivalences (if we are in a setting where this makes sense, such as when working with cochain complexes of sheaves on $X_\bullet$, as in \cref{section:cosimplicial-diagrams-of-categories}), and denote the (sub)category of such objects by $\cartSh(X_\bullet)$.
\end{definition}

\begin{remark}
  In \cref{definition:naive-definition-of-sheaf-on-simplicial-space}, we require functorial morphisms $\cal{F}^\bullet\alpha$ for \emph{all} morphisms $\alpha\colon[p]\to[q]$ in $\Delta$, but in \cite{Gre1980} these are only required to be given for all \emph{coface} morphisms $\alpha\colon[p]\to[q]$ in $\Delta_+$.
  This is entirely analogous to how a simplicial space $X_\bullet$ has a geometric realisation $|X_\bullet|$, given by a certain coend/quotient over all morphisms in $\Delta$, but also a \emph{fat} geometric realisation $\|X_\bullet\|$, where the coend/quotient only takes coface maps into consideration (and is thus defined for mere semi-simplicial spaces).
  In good scenarios (i.e.  when $X_\bullet$ is a \emph{good} simplicial space), the natural morphism $\|X_\bullet\|\to|X_\bullet|$ is a homotopy equivalence.
  One can imagine repeating all the statements in this current paper for semi-simplicial sheaves, sheaves on semi-simplicial spaces, etc., and obtaining some similar sort of equivalence result.
  We do not, however, concern ourselves with that here, although we do work with such semi-simplicial objects (for mere convenience) in \cref{subsection:certain-string-diagrams}.
\end{remark}

At a first glance, it seems like we cannot simply say that a sheaf on $X_\bullet$ is a (co)simplicial object in some category of sheaves, since each $\mathscr{F}^p$ is a sheaf on a different space (namely $X_p$).
However, thanks to the Grothendieck construction we actually \emph{can} write such objects as simplicial objects in a single category.
We follow \cite{Spi2020} in calling the resulting category the \emph{lens category}.

\begin{definition}
  Let $F\colon\cal{C}^\op\to\Cat$ be a functor.
  We define the category $\Lens{F}$ by the following data:
  \begin{itemize}
    \item an object is a pair $(c,x)$, where $c\in\cal{C}$ and $x\in F(c)$;
    \item a morphism $f\colon(c,x)\to(c',x')$ is a pair $f=(f_0,f^\sharp)$, where $f_0\colon c\to c'$ in $\cal{C}$ and $f^\sharp\colon (Ff_0)(x')\to x$ in $F(c)$.
  \end{itemize}
  We denote by $\pi\colon\Lens{F}\to\cal{C}$ the projection functor given by $\pi(c,x) := c$.
\end{definition}

The prototypical lens category for us\footnote{For others, it might be the category of \emph{polynomial functors}, obtained by taking $F$ to be the functor $\Set^\op\to\Cat$ defined by $X\mapsto\Set/X$.} is given by the following example.

\begin{example}
  The Grothendieck construction applied to the functor
  \[
    \begin{aligned}
      \Sh\colon \Space^\op &\longrightarrow \Cat
    \\X &\longmapsto \Sh(X)
    \\(f\colon X\to Y) &\longmapsto (f^*\colon\Sh(Y)\to\Sh(X))
    \end{aligned}
  \]
  gives us the category of ``all (structure) sheaves on all spaces'':
  an object in $\Lens{\Sh}$ is a pair $(X,\cal{F})$, where $X\in\Space$ and $\cal{F}\in\Sh(X)$;
  a morphism $f\colon(X,\cal{F})\to(Y,\cal{G})$ is a pair $f=(f_0,f^\sharp)$, where $f_0\colon X\to Y$ is a continuous function and $f^\sharp\colon f_0^*\cal{G}\to\cal{F}$ is a morphism of sheaves on $X$.
\end{example}

We could further refine the above example, for instance, to the case of sheaves of modules on locally ringed spaces, simply by considering the functor $\Sh_\mathrm{ring}$ that sends a ringed space $(X,\cal{O}_X)$ to the category of sheaves of $\cal{O}_X$-modules on $X$ (cf. \cite[Example~3.8]{Spi2020}).

The reason that we are so interested in this specific construction is the following lemma, which follows directly from writing out the definitions.

\begin{lemma}
  The objects of the category $\simplicial{\Lens{\Sh}}$ of simplicial objects in $\Lens{\Sh}$ are exactly sheaves on simplicial spaces.
  A sheaf $\cal{F}^\bullet\in\Sh(X_\bullet)$ on a given simplicial space $X_\bullet$ is exactly a lift
  \[
    \begin{tikzcd}
      & \Lens{\Sh} \dar["\pi"]
    \\\Delta^\op \rar[swap,"X_\bullet"] \urar[dashed,"\cal{F}^\bullet"]
      & \Space
    \end{tikzcd}
  \]
  We can further recover the notion of cartesian (resp. strictly cartesian) sheaves by simply applying the lens construction to the refined functor that sends $X$ to category of sheaves on $X$ with morphisms given only by the weak equivalences (resp. only by the isomorphisms).
\end{lemma}

\begin{remark}
  In terms of \emph{objects}, we can think of $\simplicial{\Lens{\Sh}}$ as the union of the categories $\Sh(X_\bullet)$ for all simplicial spaces $X_\bullet$, but there is a difference in the \emph{morphisms}: a morphism $\cal{F}^\bullet\to\cal{G}^\bullet$ in $\Sh(X_\bullet)$ consists of morphisms $\varphi^p\colon\cal{F}^p\to\cal{G}^p$ of sheaves on $X_p$; a morphism $(X_\bullet,\cal{F}^\bullet)\to(Y_\bullet,\cal{G}^\bullet)$ in $\simplicial{\Lens{\Sh}}$ consists of a morphism $f_\bullet\colon X_\bullet\to Y_\bullet$ of spaces along with morphisms $\varphi^p\colon f^*\cal{G}^p\to\cal{F}^p$ of sheaves on $X_p$.
  In particular, to recover $\Hom_{\Sh(X_\bullet)}(\cal{F}^\bullet,\cal{G}^\bullet)$ we need to take the \emph{opposite} of $\Hom_{\simplicial{\Lens{\Sh}}}((X_\bullet,\cal{F}^\bullet),(X_\bullet,\cal{G}^\bullet))$.

  But this is nothing specific to (co)simplicial sheaves --- this difference in direction of morphisms is already present when looking at the category of sheaves on a fixed space, where $\Hom(\cal{F},\cal{G})$ consists of maps $\varphi\colon\cal{F}\to\cal{G}$, and the category of ringed spaces, where $\Hom((X,\cal{O}_X),(Y,\cal{O}_Y))$ consists of maps $f\colon X\to Y$ along with maps $\varphi\colon f^*\cal{G}\to\cal{F}$.
\end{remark}

\subsection{The relation between the two}

As mentioned in \cref{remark:sheaves-on-simplicial-spaces-look-cosimplicial}, sheaves on simplicial spaces look somewhat like cosimplicial objects.
Indeed, it turns out that the sequence ``sheaf, \emph{co}simplicial sheaf, sheaf on a simplicial space'' consists of three constructions, each of which strictly generalises the one before.
Because of this, we will only really consider \emph{co}simplicial sheaves from now on; we do not feel negligent in not discussing \emph{simplicial} (pre)sheaves, since there is a plethora of good references on these objects (and their homotopy theory) already (e.g.  \cite{Jar2001}).

Note, first of all, that we can write the category of cosimplicial sheaves using the Grothendieck construction, as we did for sheaves on simplicial spaces.
Indeed, consider the functor
\[
  \begin{aligned}
    \cSh\colon \Space^\op &\longrightarrow \Cat
  \\X &\longmapsto \cSh(X)
  \\(f\colon X\to Y) &\longmapsto (f^*\colon\cSh(Y)\to\cSh(X)).
  \end{aligned}
\]
Then, again as an immediate consequence of the definitions, we have the following lemma.

\begin{lemma}
  The objects of the category $\Lens{\cSh}$ are exactly the cosimplicial sheaves on spaces.
  A cosimplicial sheaf $\cal{F}^\bullet\in\cSh(X)$ on a given space $X$ is exactly a lift
  \[
    \begin{tikzcd}
      & \Lens{\cSh} \dar["\pi"]
    \\\{*\} \rar[swap,"X"] \urar[dashed,"\cal{F}^\bullet"]
      & \Space.
    \end{tikzcd}
  \]
\end{lemma}

So we have two lens categories in particular:
\begin{enumerate}
  \item the category $\Lens{\cSh}$ of cosimplicial sheaves on spaces, whose objects are exactly the objects of $\cSh(X)$ for all spaces $X$; and
  \item the category $\simplicial{\Lens{\Sh}}$ of sheaves on simplicial spaces, whose objects are exactly the objects of $\Sh(X_\bullet)$ for all simplicial spaces $X_\bullet$.
\end{enumerate}
We now describe how they are related to one another, starting with the fixed space case, i.e.  $\cSh(X)$ and $\Sh(X_\bullet)$, before moving on to the Grothendieck constructions, i.e.  $\Lens{\cSh}$ and $\simplicial{\Lens{\Sh}}$.

\begin{lemma}
\label{lemma:cosimplicial-sheaf-is-sheaf-on-constant}
  A cosimplicial sheaf on a space is exactly a sheaf on the corresponding constant simplicial space.
  That is, given a space $X$, we have an equivalence of categories
  \[
    \cSh(X) \simeq \Sh(\constant{X}_\bullet)
  \]
  where $\constant{X}_\bullet$ is the simplicial space with $\constant{X}_p=X$ for all $p\in\mathbb{N}$ and $\constant{X}_\bullet(\alpha)=\id_X$ for all $\alpha\colon[p]\to[q]$ in $\Delta$.
\end{lemma}

\begin{proof}
  Let $\cal{F}^\bullet\in\Sh(\constant{X}_\bullet)$.
  Then we have $\cal{F}^p\in\Sh(\constant{X}_p)$ for all $p\in\mathbb{N}$, along with functorial $\cal{F}^\bullet\alpha\colon(X_\bullet\alpha)^*\cal{F}^p\to\cal{F}^q$ for all $\alpha\colon[p]\to[q]$ in $\Delta$.
  But, by definition, $\constant{X}_p=X$ for all $p\in\mathbb{N}$, and $\constant{X}_\bullet\alpha=\id_X$ for all $\alpha$, so $\cal{F}^\bullet\alpha$ is really a morphism $\cal{F}^p\to\cal{F}^q$ in $\Sh(X)$.
  That is, objects of $\Sh(\constant{X}_\bullet)$ are exactly objects of $\cSh(X)$.
  By the same arguments, the morphisms in the two categories also agree on the nose.
\end{proof}

\begin{lemma}
\label{lemma:const-gives-reflective-subcategory}
  The functor $\const\colon\Space\to\sSpace$ given by $X\mapsto\constant{X}_\bullet$ induces a fully faithful functor
  \[
    \begin{aligned}
      \const\colon\Lens{\cSh} &\longhookrightarrow \simplicial{\Lens{\Sh}}
    \\(X,\cal{F}^\bullet) &\longmapsto (\constant{X}_\bullet,\cal{F}^\bullet).
    \end{aligned}
  \]
  Further, this inclusion $\const$ has a left adjoint
  \[
    \begin{aligned}
      \colim\colon\simplicial{\Lens{\Sh}} &\longrightarrow \Lens{\cSh}
    \\(X_\bullet,\cal{F}^\bullet) &\longmapsto (\colim X_\bullet, \{(X_p\to\colim X_\bullet)_*\cal{F}^p\}_{p\in\mathbb{N}})
    \end{aligned}
  \]
  given by taking the colimit of the simplicial space and the pushforwards of the sheaves along the inclusions into the colimit, thus witnessing $\Lens{\cSh}$ as a reflective subcategory of $\simplicial{\Lens{\Sh}}$.
\end{lemma}

\begin{proof}
  This is a reasonably immediate consequence of the fact that pullback is left adjoint to pushforward, combined with the universal property of the colimit, but we spell out the details here just to be clear.
  Consider a morphism
  \[
    (f_\bullet,\varphi^\bullet)
    \in \Hom_{\simplicial{\Lens{\Sh}}}
    \big(
      (X_\bullet,\cal{F}^\bullet),(\constant{Y}_\bullet,\cal{G}^\bullet)
    \big)
  \]
  so that $f_\bullet\colon X_\bullet\longrightarrow\constant{Y}_\bullet$ and $\varphi^\bullet\colon f_\bullet^*\cal{G}^\bullet\longrightarrow\cal{F}^\bullet$.
  By the universal property of the colimit, this induces a unique morphism
  \[
    |f_\bullet|\colon\colim X_\bullet\longrightarrow Y
  \]
  along with morphisms
  \[
    i_p\colon X_p\longrightarrow\colim X_\bullet
  \]
  such that $f_p=|f_\bullet|\circ i_p$.
  But pullback is left adjoint to pushforward, and so the data of a morphism
  \[
    \varphi^p\colon f_p^*\cal{G}^p\longrightarrow\cal{F}^p
  \]
  is equivalent to the data of a morphism
  \[
    \psi_p\colon\cal{G}^p\longrightarrow (f_p)_*\cal{F}^p=(|f_\bullet|)_*(i_p)_*\cal{F}^p
  \]
  and thus to the data of a morphism
  \[
    \widetilde{\varphi}^p\colon|f_\bullet|^*\cal{G}^p\longrightarrow(i_p)_*\cal{F}^p
  \]
  which, all together, gives
  \[
    (|f_\bullet|,\widetilde{\varphi}^p)
    \in \Hom_{\Lens{\cSh}}
    \big(
      (\colim X_\bullet,(i_\bullet)_*\cal{F}^\bullet),(Y,\cal{G}^\bullet)
    \big)
  \]
  which witnesses the necessary (natural) isomorphism for the adjunction.
\end{proof}

\begin{remark}
  There is a dual statement of \cref{lemma:const-gives-reflective-subcategory}, where the limit (along with the pushforwards along the maps from the limit) gives a right adjoint to the constant functor, thus witnessing $\Lens{\cSh}$ as a \emph{co}reflective subcategory of $\simplicial{\Lens{\Sh}}$ as well.
  We do not, however, care so much about this limit construction here, since the main case of interest for us is where $X_\bullet=X_\bullet^\cover$ is the Čech nerve of a space, and the geometric realisation is given by the \emph{co}limit, i.e.  $|X_\bullet^\cover| :=\colim X_\bullet^\cover\simeq X$.
\end{remark}

Analogously to how cosimplicial sheaves are specific examples of sheaves on a simplicial space, we can actually recover sheaves as specific examples of cosimplicial sheaves, as follows.

\begin{lemma}
\label{lemma:sheaf-is-strictly-cartesian-sheaf-on-constant}
  A sheaf on a space is exactly a strictly cartesian sheaf on a constant simplicial space.
  That is, given a space $X$, we have an equivalence of categories
  \[
    \Sh(X) \simeq \scartSh(\constant{X}_\bullet)
  \]
  where $\constant{X}_\bullet$ is the constant simplicial space, as in \cref{lemma:cosimplicial-sheaf-is-sheaf-on-constant}.
\end{lemma}

\begin{proof}
  We already know by \cref{lemma:cosimplicial-sheaf-is-sheaf-on-constant} that an object $\cal{F}^\bullet\in\Sh(\constant{X}_\bullet)$ is exactly an object of $\cSh(X)$.
  By definition, $\cal{F}^\bullet$ being strictly cartesian means that the $\cal{F}^\bullet\alpha\colon\cal{F}^p\to\cal{F}^q$ are isomorphisms for all $\alpha\colon[p]\to[q]$ in $\Delta$, and so the functor {$\scartSh(\constant{X}_\bullet)\to\Sh(X)$} given by $\cal{F}^\bullet\mapsto\cal{F}^0$ has an up-to-isomorphism inverse given by $\cal{F}^0\mapsto(\constant{\cal{F}^0})^\bullet$, where $(\constant{\cal{F}^0})^p :=\cal{F}^0$ for all $p\in\mathbb{N}$, and $(\constant{\cal{F}^0})^\bullet\alpha :=\id_{\cal{F}^0}$ for all $\alpha\colon[p]\to[q]$ in $\Delta$.
\end{proof}

\begin{table}[h!]
    \centering
    \begin{tabular}{clc}
      Lens category
      & Objects
      & Single fixed space
    \\\midrule
      $\Lens{\Sh}$
      & sheaves on spaces
      & $\Sh(X)\simeq\scartSh(\constant{X}_\bullet)$
    \\$\Lens{\cSh}$
      & cosimplicial sheaves on spaces
      & $\cSh(X)\simeq\Sh(\constant{X}_\bullet)$
    \\$\simplicial{\Lens{\Sh}}$
      & sheaves on simplicial spaces
      & $\Sh(X_\bullet)$
    \end{tabular}
    \caption{A summary of the three notions studied so far, with the equivalences in the last column following from \cref{lemma:sheaf-is-strictly-cartesian-sheaf-on-constant,lemma:cosimplicial-sheaf-is-sheaf-on-constant}.}
    \label{table:summary-of-three-notions}
\end{table}

\begin{figure}[h!]
    \centering
    \begin{tikzpicture}[scale=0.9]
      \tikzstyle{space}=[inner sep=0.25cm,draw,thick]
      \tikzstyle{ellipsish}=[rounded corners=0.35cm]
      \tikzstyle{shorter}=[shorten >=0.2cm,shorten <=0.2cm]
      % FIGURE 1: sheaf on a space
      \begin{scope}
        % label
        \node at (0,-1.5) {$\cal{F}\in\Sh(X)$};
        % sheaves and bases
        \node[space,ellipsish] (sheaf) at (0,2) {$\cal{F}$};
        \node[space] (base) at (0,0) {$X$};
        % arrows from sheaves to bases
        \draw[thick,-latex,shorter] (sheaf) to (base);
      \end{scope}
      % FIGURE 2: cosimplicial sheaf on a space
      \begin{scope}[shift={(5,0)}]
        % label
        \node at (0,-1.5) {$\cal{F}^\bullet\in\cSh(X)$};
        % sheaves and bases
        \node[space,ellipsish] (sheaf0) at (-2,2) {$\cal{F}^0$};
        \node[space,ellipsish] (sheaf1) at (0,2) {$\cal{F}^1$};
        \node (sheaf2) at (2,2) {\ldots};
        \node[space] (base) at (0,0) {$X$};
        % arrows from sheaves to bases
        \draw[thick,-latex,shorter,bend right] (sheaf0) to (base.west);
        \draw[thick,-latex,shorter] (sheaf1) to (base);
        \draw[thick,-latex,shorter,bend left] (sheaf2) to (base.east);
        % cosimplicial maps upstairs from 0 to 1
        \draw[thick,->,shorter,transform canvas={yshift=2mm}] (sheaf0) to (sheaf1);
        \draw[thick,->,shorter] (sheaf1) to (sheaf0);
        \draw[thick,->,shorter,transform canvas={yshift=-2mm}] (sheaf0) to (sheaf1);
        % cosimplicial maps upstairs from 1 to 2
        \draw[thick,->,shorter,transform canvas={yshift=4mm}] (sheaf1) to (sheaf2);
        \draw[thick,->,shorter,transform canvas={yshift=2mm}] (sheaf2) to (sheaf1);
        \draw[thick,->,shorter] (sheaf1) to (sheaf2);
        \draw[thick,->,shorter,transform canvas={yshift=-2mm}] (sheaf2) to (sheaf1);
        \draw[thick,->,shorter,transform canvas={yshift=-4mm}] (sheaf1) to (sheaf2);
      \end{scope}
      % FIGURE 3: sheaf on a simplicial space
      \begin{scope}[shift={(11.5,0)}]
        % label
        \node at (0,-1.5) {$\cal{F}^\bullet\in\Sh(X_\bullet)$};
        % sheaves and bases
        \node[space,ellipsish] (sheaf0) at (-2,2) {$\cal{F}^0$};
        \node[space,ellipsish] (sheaf1) at (0,2) {$\cal{F}^1$};
        \node (sheaf2) at (2,2) {\ldots};
        \node[space] (base0) at (-2,0) {$X_0$};
        \node[space] (base1) at (0,0) {$X_1$};
        \node (base2) at (2,0) {\ldots};
        % arrows from sheaves to bases
        \draw[thick,-latex,shorter] (sheaf0) to (base0);
        \draw[thick,-latex,shorter] (sheaf1) to (base1);
        \draw[thick,-latex,shorter] (sheaf2) to (base2);
        % cosimplicial maps upstairs from 0 to 1
        \draw[thick,->,shorter,transform canvas={yshift=2mm}] (sheaf0) to (sheaf1);
        \draw[thick,->,shorter] (sheaf1) to (sheaf0);
        \draw[thick,->,shorter,transform canvas={yshift=-2mm}] (sheaf0) to (sheaf1);
        % cosimplicial maps upstairs from 1 to 2
        \draw[thick,->,shorter,transform canvas={yshift=4mm}] (sheaf1) to (sheaf2);
        \draw[thick,->,shorter,transform canvas={yshift=2mm}] (sheaf2) to (sheaf1);
        \draw[thick,->,shorter] (sheaf1) to (sheaf2);
        \draw[thick,->,shorter,transform canvas={yshift=-2mm}] (sheaf2) to (sheaf1);
        \draw[thick,->,shorter,transform canvas={yshift=-4mm}] (sheaf1) to (sheaf2);
        % simplicial maps downstairs from 0 to 1
        \draw[thick,->,shorter,transform canvas={yshift=2mm}] (base1) to (base0);
        \draw[thick,->,shorter] (base0) to (base1);
        \draw[thick,->,shorter,transform canvas={yshift=-2mm}] (base1) to (base0);
        % simplicial maps downstairs from 1 to 2
        \draw[thick,->,shorter,transform canvas={yshift=4mm}] (base2) to (base1);
        \draw[thick,->,shorter,transform canvas={yshift=2mm}] (base1) to (base2);
        \draw[thick,->,shorter] (base2) to (base1);
        \draw[thick,->,shorter,transform canvas={yshift=-2mm}] (base1) to (base2);
        \draw[thick,->,shorter,transform canvas={yshift=-4mm}] (base2) to (base1);
      \end{scope}
    \end{tikzpicture}
    \caption{A pictorial representation of how sheaves, cosimplicial sheaves, and sheaves on simplicial spaces relate to one another, and are of strictly increasing generality. Note that the ``squares'' in the diagram on the far right (for $\cal{F}^\bullet\in\Sh(X_\bullet)$) do \emph{not} commute \emph{unless} the sheaf is strongly cartesian. This is the first hint at the fact that these sheaves arise as \emph{lax} limit objects, as we will formalise in \cref{lemma:laxholim-and-holim-of-sh}.}
    \label{figure:picture-of-three-notions}
\end{figure}
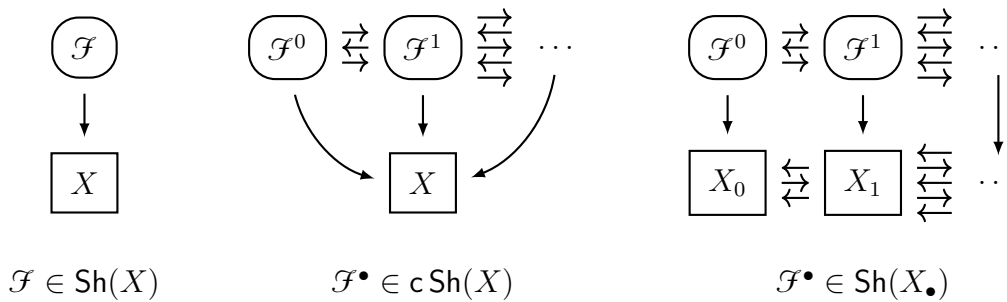

\begin{remark}
\label{remark:the-four-notions}
  In both \cref{table:summary-of-three-notions} and \cref{figure:picture-of-three-notions}, we are omitting a fourth family of objects, namely strictly cartesian sheaves on \emph{not-necessarily constant} simplicial spaces: $\scartSh(X_\bullet)$.
  If we were to draw these as in \cref{figure:picture-of-three-notions}, they would look as follows:
  
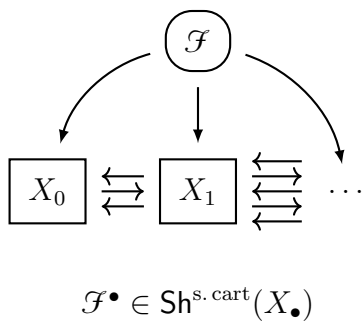
\[
    \begin{tikzpicture}
      \tikzstyle{space}=[inner sep=0.25cm,draw,thick]
      \tikzstyle{ellipsish}=[rounded corners=0.35cm]
      \tikzstyle{shorter}=[shorten >=0.2cm,shorten <=0.2cm]
      % FIGURE 4: strictly carsetian sheaf on a simplicial space
      \begin{scope}[shift={(11.5,0)}]
        % label
        \node at (0,-1.5) {$\cal{F}^\bullet\in\scartSh(X_\bullet)$};
        % sheaves and bases
        \node[space,ellipsish] (sheaf1) at (0,2) {$\cal{F}$};
        \node[space] (base0) at (-2,0) {$X_0$};
        \node[space] (base1) at (0,0) {$X_1$};
        \node (base2) at (2,0) {\ldots};
        % arrows from sheaves to bases
        \draw[thick,-latex,shorter,bend right] (sheaf1) to (base0);
        \draw[thick,-latex,shorter] (sheaf1) to (base1);
        \draw[thick,-latex,shorter,bend left] (sheaf1) to (base2);
        % simplicial maps downstairs from 0 to 1
        \draw[thick,->,shorter,transform canvas={yshift=2mm}] (base1) to (base0);
        \draw[thick,->,shorter] (base0) to (base1);
        \draw[thick,->,shorter,transform canvas={yshift=-2mm}] (base1) to (base0);
        % simplicial maps downstairs from 1 to 2
        \draw[thick,->,shorter,transform canvas={yshift=4mm}] (base2) to (base1);
        \draw[thick,->,shorter,transform canvas={yshift=2mm}] (base1) to (base2);
        \draw[thick,->,shorter] (base2) to (base1);
        \draw[thick,->,shorter,transform canvas={yshift=-2mm}] (base1) to (base2);
        \draw[thick,->,shorter,transform canvas={yshift=-4mm}] (base2) to (base1);
      \end{scope}
    \end{tikzpicture}
  \]
  Although these objects are somewhat too strict to be of immediate interest, the \emph{weakly cartesian} version proves to be incredibly useful, as explained in the context of complex geometry in \cref{subsection:complex-geometry}.
\end{remark}

\begin{remark}
  In the case where $X_\bullet=X_\bullet^\cover$ is the Čech nerve of a cover, all four notions fit into the diagram
  \[
    \begin{tikzcd}
      &\Sh(X)
        \ar[dl,hook']
        \ar[dr,hook,two heads]
    \\\cSh(X)
        \ar[dr,hook]
      &&\scartSh(X_\bullet^\cover)
        \ar[dl,hook']
    \\&\Sh(X_\bullet^\cover)
    \end{tikzcd}
  \]
  where the top-right arrow is an equivalence (since the functor $\Sh\colon\Space^\op\to\Cat$ is a stack), and the bottom-right arrow is the natural map from the lax (homotopy) limit into the (homotopy) limit (cf. \cref{lemma:laxholim-and-holim-of-sh}).
  Note that the parallel arrows are ``the same'', i.e.  the two arrows that go down and to the right are given by ``replace $X$ with $X_\bullet^\cover$'', and the two arrows that go down and to the left are given by ``take cosimplicial objects in the sheaf part''.
\end{remark}

\section{Sections of sheaves}
\label{section:sections-of-sheaves}

Although simplicial sheaves (in all meanings of the phrase) have been discussed elsewhere in the literature, we are not familiar with any study of what a \emph{section} of such a sheaf should be.
Here we propose a definition.
Note, however, that simplicial sheaves (in all meanings of the phrase) should form a topos, and from this we could recover a definition of section using the general theory of topoi.

\subsection{Sections of cosimplicial sheaves}

We start by giving a definition, and then justify why it is indeed ``the good one''.

\begin{definition}
\label{definition:section-of-cosimplicial-sheaf}
  Let $\cal{F}^\bullet\in\cSh(X)$ be a cosimplicial sheaf on a space $X$.
  Then we define the \emph{global sections} of $\cal{F}^\bullet$ to be
  \[
    \cal{F}^\bullet(X)
    =
    \Gamma(X,\cal{F}^\bullet)
     :=
    \eq\left(
      \cal{F}^0(X)\rightrightarrows\cal{F}^1(X).
    \right)
  \]
  i.e.  a global section of $\cal{F}$ is an element of the equaliser of the $1$-truncation of $\cal{F}^\bullet$.
\end{definition}

Forgetting the above for a moment, the definition that might appear to be the most natural would be that which follows from the identification of $\cosimplicial{\PSh(X)}$ with $[X^\op,\cSet]$, i.e.\
\begin{quotation}
    \noindent
    A global section of $\cal{F}^\bullet\in\cSh(X)$ is a collection of sections $(\sigma^p\in\Gamma(X,\cal{F}^p))_{p\in\mathbb{N}}$ such that $\cal{F}^\bullet\alpha(\sigma^p)=\sigma^q$ for all $\alpha\colon[p]\to[q]$ in $\Delta$.
\end{quotation}
But since this is a strict condition (i.e.  we are asking for \emph{equality}), each $\sigma^p$ is determined by $\{\sigma^0,\ldots,\sigma^{p-1}\}$, so we need only ensure the initial data of $\sigma^0$ is such that we can construct $\sigma^1$, i.e.  a section should be exactly some $\sigma^0\in\eq(\cal{F}^0(X)\rightrightarrows\cal{F}^1(X))$, and this recovers \cref{definition:section-of-cosimplicial-sheaf}.

\begin{definition}
  Let $\cal{F}^\bullet\in\cSh(X)$ be a cosimplicial sheaf on a space $X$, and $i\colon U\hookrightarrow X$ the inclusion of an open subset.
  Then we define the \emph{local sections} of $\cal{F}^\bullet$ (over $U$) to be the global sections of $i^*\cal{F}^\bullet\in\cSh(U)$.
\end{definition}

\subsection{Sections of sheaves on a simplicial space}

The definition for arbitrary sheaves on a simplicial space is morally the same as the one for cosimplicial sheaves given in \cref{definition:section-of-cosimplicial-sheaf} (following the fact that cosimplicial sheaves are ``just'' sheaves on a constant simplicial space, by \cref{lemma:cosimplicial-sheaf-is-sheaf-on-constant}), but since $X_\bullet$ is no longer a priori constant, $\cal{F}^0$ and $\cal{F}^1$ live on different spaces, so we have to take pullbacks along the two coface maps.

\begin{definition}
\label{definition:section-of-sheaf-on-simplicial-space}
  Let $\cal{F}^\bullet\in\Sh(X_\bullet)$ be a sheaf on a simplicial space $X_\bullet$.
  Then we define the \emph{global sections} of $\cal{F}^\bullet$ to be
  \[
    \cal{F}^\bullet(X_\bullet)
    =
    \Gamma(X_\bullet,\cal{F}^\bullet)
     :=
    \lim\left(
      \begin{tikzcd}[row sep=tiny,column sep=5em]
        \left(\left(X_\bullet f_1^0\right)^*\cal{F}^0\right)(X_1)
          \ar[dr,near start,"\left(\cal{F}^\bullet f_1^0\right)(X_1)"]
      \\&\cal{F}^1(X_1)
      \\\left(\left(X_\bullet f_1^1\right)^*\cal{F}^0\right)(X_1)
          \ar[ur,near start,swap,"(\cal{F}^\bullet f_1^1)(X_1)"]
      \end{tikzcd}
    \right)
  \]
  where $f_1^i\colon[0]\to[1]$ is the $i$th coface map in $\Delta$.
\end{definition}

\begin{definition}
  Let $\cal{F}^\bullet\in\Sh(X_\bullet)$ be a sheaf on a simplicial space $X_\bullet$, and \mbox{$i_\bullet\colon U_\bullet\hookrightarrow X_\bullet$} the inclusion of a simplicial space such that each $i_p\colon U_p\hookrightarrow X_p$ is the inclusion of an open subset.
  Then we define the \emph{local sections} of $\cal{F}^\bullet$ (over $U_\bullet$) to be the global sections of $i_\bullet^*\cal{F}^\bullet\in\Sh(U_\bullet)$.
\end{definition}

\subsection{Weak sections}
\label{subsection:weak-sections}

Recall that, given a sheaf $\cal{F}$ (of sets) on a space $X$, its global sections $\cal{F}(X)$ are in bijection with morphisms in the $\Hom$-set from the constant sheaf on a singleton set, i.e.\
\[
  \cal{F}(X) = \Gamma(X,\cal{F}) \cong \Hom_{\Sh(X)}(\constant{\{*\}},\cal{F})
\]
(which is also equivalent to the direct image along the morphism $X\to\{*\}$).
It turns out that \cref{definition:section-of-cosimplicial-sheaf,definition:section-of-sheaf-on-simplicial-space} can also be written in this form: as functors from $\cal{C}$ to $\Set$, where $\cal{C}$ is either $\cSh(X)$ or $\Sh(X_\bullet)$, we have an isomorphism
\[
  \Gamma(X,-) \cong \Hom_{\cal{C}}(\constant{\{*\}},-).
\]
But, as we mentioned, with this definition a section $\sigma^\bullet$ ends up being determined entirely by its degree-$0$ part, since we have conditions of strict equality.
If we wish to recover a weaker notion of section (in the case of $\cSh(X)$, say), and we were working with sheaves of objects that have some notion of weak equivalence, then we could take the (homotopy) limit of the \emph{entire} co(semi-)simplicial diagram
\[
  \begin{tikzcd}
    \cal{F}^0(X)
      \ar[r,shift left=1]
      \ar[r,shift right=1]
    & \cal{F}^1(X)
      \ar[r,shift left=2]
      \ar[r]
      \ar[r,shift right=2]
    & \ldots
  \end{tikzcd}
\]
instead of just the $1$-truncation, which would correspond to replacing the constant sheaf on a singleton set $\constant{\{*\}}$ with the tower of constant sheaves on the geometric simplices $\{*\} \cong \Delta^0 \hookrightarrow \Delta^1 \hookrightarrow \cdots$. That is, we would have some section $\sigma^0\in\Gamma(X,\cal{F}^0)$ as before, but then \emph{two} (weakly equivalent) sections $\sigma_0^1\simeq\sigma_1^1\in\Gamma(X,\cal{F}^1)$, and such that $\sigma^0$ is mapped to $\sigma_i^1$ by the $i$th coface map $f_1^i\colon[0]\to[1]$, and then \emph{three} (weakly equivalent) sections $\sigma_0^2\simeq\sigma_1^2\simeq\sigma_2^2\in\Gamma(X,\cal{F}^2)$ such that \ldots, and so on.
In fact, we could even consider a \emph{lax} version of this, where we do not ask for the morphism $\sigma_0^1\to\sigma_1^1$ to be a weak equivalence, but instead simply to be an arbitrary morphism (and, indeed, we consider such a notion in \cref{subsection:certain-string-diagrams}).

This approach (via the $\Hom$ out of a point, or out of the ``fat point'' $\Delta^0\hookrightarrow\Delta^1\hookrightarrow\ldots$) also explains why sections ``look like'' elements of the totalisation of a cosimplicial simplicial set: this is a way of computing the homotopy limit of the cosimplicial simplicial model of the mapping space.

\section{Cosimplicial diagrams of categories}
\label{section:cosimplicial-diagrams-of-categories}

In this section we will consider the specific case of locally ringed spaces and sheaves of modules.
If $(X,\cal{O}_X)$ is a locally ringed space, then we simply write $\Sh(X)$ to mean \emph{the category of cochain complexes of sheaves of $\cal{O}_X$-modules} on $X$.
This is the setting of interest to most algebraic geometers, and we consider it here because it allows us to endow $\Sh(X)$ with the structure of a model category.
Another simplification we make is to only consider the case where $X_\bullet=X_\bullet^\cover$ is the Čech nerve of a locally finite cover $\cover$ of a locally ringed space $(X,\cal{O}_X)$, giving us the simplicial locally ringed space $(X_\bullet^\cover,\cal{O}_{X_\bullet^\cover})$, where $\cal{O}_{X_p^\cover} :=(\pi_p\colon X_p^\cover\to X)^{-1}\cal{O}_X$.
For clarity, we now spell out the details of \cref{definition:naive-definition-of-sheaf-on-simplicial-space} in this specific setting.

A \emph{cochain complex of sheaves on the simplicial locally ringed space $(X_\bullet^\cover,\cal{O}_{X_\bullet^\cover})$} is the data of sheaves $(\cal{F}^{p,i})_{i\in\mathbb{Z}}$ of $\cal{O}_{X_p^\cover}$-modules on $X_p^\cover$ for each $p\in\mathbb{N}$ along with functorial morphisms
\[
  \cal{F}^{\bullet,i}\alpha\colon (X_\bullet^\cover\alpha)^*\cal{F}^{p,i} \longrightarrow \cal{F}^{q,i}
\]
of $\cal{O}_{X_q^\cover}$-modules for all $\alpha\colon[p]\to[q]$ in $\Delta$, as well as morphisms
\[
  d^i =
  \left(
    d^{p,i}\colon\cal{F}^{p,i} \longrightarrow \cal{F}^{p,i+1}
  \right)_{p\in\mathbb{N}}
\]
of sheaves on $X_\bullet^\cover$ for all $i\in\mathbb{Z}$ such that
\[
  \begin{tikzcd}[sep=huge]
    (X_\bullet^\cover\alpha)^*\cal{F}^{p,i}
      \ar[r,"(X_\bullet^\cover\alpha)^* d^{p,i}"]
      \ar[d,swap,"\cal{F}^{\bullet,i}\alpha"]
    & (X_\bullet^\cover\alpha)^*\cal{F}^{p,i+1}
      \ar[d,"\cal{F}^{\bullet,i+1}\alpha"]
  \\\cal{F}^{q,i}
      \ar[r,swap,"d^{q,i}"]
    & \cal{F}^{q,i+1}
  \end{tikzcd}
\]
commutes for all $\alpha\colon[p]\to[q]$ in $\Delta$, and such that $d^{i+1}\circ d^i=0$ for all $i\in\mathbb{Z}$.
We say that such an object $\cal{F}^{\bullet,\anotherbullet}$ is (\emph{weakly}) \emph{cartesian} if the morphisms
\[
  \cal{F}^{\bullet,\anotherbullet}\alpha\colon(X_\bullet^\cover\alpha)^*\cal{F}^{p,\anotherbullet}\to\cal{F}^{q,\anotherbullet}
\]
of cochain complexes are quasi-isomorphisms\footnote{Note that the condition that the $d^i$ be morphisms of sheaves on $X_\bullet^\cover$ already ensures that the $\cal{F}^{\bullet,\anotherbullet}\alpha$ will be chain maps, i.e.  that they commute with the differentials.} for all $\alpha\colon[p]\to[q]$ in $\Delta$.

\subsection{(Lax) homotopy limits of diagrams of model categories}

\begin{definition}
  Let $\cal{D}$ be a small category
  Then a \emph{$\cal{D}$-shaped diagram of model categories} is a functor $M\colon\cal{D}\to\ModL$, where $\ModL$ is the (large) category whose objects are model categories and whose morphisms are left Quillen functors.
  That is, $M$ is the data of a model category $M_d$ for all $d\in\cal{D}$ along with a functorial\footnote{Meaning $F_{d,f}^{\eta\theta}=F_{e,f}^\eta\circ F_{d,e}^\theta$ for all $\theta\colon d\to e$ and $\eta\colon e\to f$ in $\cal{D}$.} assignment of left Quillen functors $F_{d,e}^\theta\colon M_d\to M_e$ for all $\theta\colon d\to e$ in $\cal{D}$.
\end{definition}

\begin{definition}[{\cite[Definition~3.1]{Ber2012}}]
\label{definition:lax-homotopy-limit-of-model-categories}
  Let $\cal{D}$ be a small category, and $M$ a $\cal{D}$-shaped diagram of model categories such that each $M_d$ is cofibrantly generated.
  Then a \emph{lax homotopy limit of $M$} is the data of the following:
  \begin{itemize}
    \item an object $m_d\in M_d$ for each $d\in\cal{D}$;
    \item a morphism $u_{d,e}^\theta\colon F_{d,e}^\theta(m_d)\to m_e$ in $M_e$ for each $\theta\colon d\to e$ in $\cal{D}$
  \end{itemize}
  such that
  \[
    u_{d,f}^{\eta\theta} = u_{e,f}^\eta\circ F_{e,f}^\eta(u_{d,e}^\theta).
  \]
  This defines a category, which we denote by $\laxholim_{d\in\cal{D}}M$.

  Given a lax homotopy limit of $M$, we define the \emph{homotopy limit} as the full subcategory $\holim_{d\in\cal{D}}M$ of $\laxholim_{d\in\cal{D}}M$ consisting of the objects such that every $u_{d,e}^\theta$ is a weak equivalence (in $M_e$).
\end{definition}

\begin{remark}
  As noted in \cite{Ber2012}, we can endow the category $\laxholim_{d\in\cal{D}}M$ with the structure of a model category (with weak equivalences and cofibrations given levelwise, appealing to the hypothesis that each $M_d$ is cofibrantly generated), but we \emph{cannot} a priori endow $\holim_{d\in\cal{D}}M$ with the structure of a model category, since the requirement that the $u_{d,e}^\theta$ be weak equivalences is \emph{not} preserved by arbitrary limits and colimits.
\end{remark}

\begin{lemma}
\label{lemma:laxholim-and-holim-of-sh}
  Consider the $\Delta$-shaped diagram of model categories
  \[
    \begin{tikzcd}
      \Sh(X_0^\cover)
        \ar[r,shift left=2]
        \ar[r,shift right=2]
      & \Sh(X_1^\cover)
        \ar[r,shift left=4]
        \ar[r]
        \ar[r,shift right=4]
        \ar[l]
      & \ldots
        \ar[l,shift left=2]
        \ar[l,shift right=2]
    \end{tikzcd}.
  \]
  Then we have equivalences
  \[
    \laxholim_{[p]\in\Delta}\Sh(X_p^\cover)
    \simeq \Sh(X_\bullet^\cover)
    \qquad\text{and}\qquad
    \holim_{[p]\in\Delta}\Sh(X_p^\cover)
    \simeq \cartSh(X_\bullet^\cover).
  \]
\end{lemma}

\begin{proof}
  This follows from taking $\cal{D}=\Delta$, $F_{\alpha,\beta}^\theta=(X_\bullet^\cover\theta)^*$, and $u_{\alpha,\beta}^\theta=\cal{F}^\bullet\theta$ in \cref{definition:lax-homotopy-limit-of-model-categories}.
\end{proof}

Note that the homotopy limit in \cref{lemma:laxholim-and-holim-of-sh} is \emph{also} equivalent to $\Sh(X)$, since the functor $\Sh\colon\Space^\op\to\Cat$ is a stack, whence there is an equivalence $\Sh(X)\simeq\cartSh(X_\bullet^\cover)$.
Indeed, the only one of the four categories in \cref{remark:the-four-notions} which is absent here is $\cSh(X)\simeq\Sh(\constant{X}_\bullet)$.
Note, however, that $X_\bullet^\cover\simeq \constant{X}_\bullet$ does \emph{not} imply that $\Sh(X_\bullet^\cover)\simeq\Sh(\constant{X}_\bullet)$, since $\Sh\colon\Space^\op\to\Cat$ is not homotopy invariant, but whether or not these two categories actually \emph{are} equivalent for some other reason is a question that we will not try to answer here.

One might wonder what, after all, the use of $\cartSh(X_\bullet^\cover)$ is, if it happens to simply be equivalent to the (much simpler) $\Sh(X)$.
The answer is very similar to why $X_\bullet^\cover$ is useful, even though its geometric realisation is equivalent to the much-simpler object $X$, in that decomposing an object into ``simplicial levels'' allows us to prescribe much finer homotopical data to different parts, instead of having to treat the object as one single lump thing. %too informal?
For a less hand-wavy justification, see \cref{subsection:complex-geometry}, which gives an example of how this simplicial decomposition lets us describe homotopical versions of coherent sheaves.

\subsection{Locally free sheaves and perfect complexes}

Before concluding this section, we give a brief aside on two variations of \cref{lemma:laxholim-and-holim-of-sh}.
Since we are working in the setting of sheaves of modules, we could also consider specific families of sheaves, such as \emph{locally free} sheaves, or \emph{perfect complexes}.
Rather than working with model categories, we can consider the dg-categories $\mathsf{LocFree}(X)$ of complexes of locally free sheaves on $X$, and $\mathsf{Perf}(X)$ of perfect complexes on $X$.
It turns out that the analogous homotopy limits for these two categories recover the notion of (\emph{perfect}) \emph{twisting cochains}:
\[
  \holim_{[p]\in\Delta}\mathsf{LocFree}\,(X_p^\cover)\simeq\mathsf{Tw}\,(X)
  \qquad\qquad
  \holim_{[p]\in\Delta}\mathsf{Perf}\,(X_p^\cover)\simeq\mathsf{Tw}_\mathrm{perf}\,(X)
\]
(by \cite[Corollary~3 and Proposition~11]{BHW2017}, respectively).

For an overview of twisting cochains and twisted complexes in algebraic geometry and topology, see e.g.  \cite[Appendix~G]{Hos2020}.

\section{Applications}
\label{section:applications}

\subsection{Complex geometry}
\label{subsection:complex-geometry}

Given a (locally) ringed space $(X,\cal{O}_X)$, geometers are often interested in a specific class of sheaves of $\cal{O}_X$-modules, namely \emph{coherent sheaves}.
These are sheaves that satisfy a nice finiteness property: around any point, they are surjected onto by a (finite) free module, and, on any open subset, the kernel of any surjection from a (finite) free module is \emph{also} locally surjected onto by a (finite) free module.
In the setting of algebraic geometry, given a commutative ring $R$, there is a nice correspondence between $R$-modules and sheaves on the Noetherian affine scheme $\Spec(R)$, summarised in the following dictionary:
\[
  \begin{aligned}
    \mbox{arbitrary $R$-modules}
    &\longleftrightarrow \mbox{quasi-coherent sheaves on $\Spec(R)$}
  \\\mbox{finitely generated $R$-modules}
    &\longleftrightarrow \mbox{coherent sheaves on $\Spec(R)$}
  \\\mbox{finitely generated projective $R$-modules}
    &\longleftrightarrow \mbox{locally free sheaves on $\Spec(R)$.}
  \end{aligned}
\]

There are two results concerning coherent sheaves in \emph{algebraic} geometry that interest us here specifically:
\begin{enumerate}
  \item On any (quasi-projective Noetherian) scheme $X$, a coherent sheaf can be resolved by locally free sheaves, i.e.  is quasi-isomorphic to a complex of locally free sheaves.
  \item The category $\Coh(X)$ of complexes of coherent sheaves (on some fixed $X$ as above) is equivalent to the category $\CCoh(X)$ of complexes of sheaves whose internal cohomology is degree-wise coherent, i.e.  complexes $(\cal{F}^\bullet,d^\bullet)$ such that the sheaves $\Ker d^i/\Im d^{i-1}$ are coherent for all $i\in\mathbb{Z}$.
\end{enumerate}

The first of these statements is \emph{not} true in the analytic case, i.e.  there exist coherent sheaves (over the sheaf of holomorphic functions) on complex-analytic manifolds which \emph{cannot} be resolved by locally free sheaves \cite[Corollary~A.5]{Voi2002}.
However, a beautiful construction in \cite{Gre1980} shows that coherent analytic sheaves can \emph{always} be \emph{locally} resolved by locally free sheaves; in the language of this current survey, this construction was used in \cite{Hos2020a} to show that (complexes of) coherent analytic sheaves are equivalent (in some suitable $(\infty,1)$-categorical sense) to a certain full subcategory of $\cartSh(X_\bullet^\cover)$.
That is, a nice subclass of complexes of cartesian sheaves on the Čech nerve gives a homotopical presentation of (complexes of) coherent analytic sheaves, modulo a caveat explained in the next paragraph.

The second of these statements is an open problem in the analytic case.
It is known that $\Coh(X)$ and $\CCoh(X)$ are equivalent in the specific case where $X$ is a smooth complex-analytic \emph{surface} (cf. \cite[Corollary~5.2.2]{BV2003}), but a general result for arbitrary (or even some fixed higher) dimensions is not known.
The specific category used in \cite{Hos2020a} lies a priori between $\Coh(X)$ and $\CCoh(X)$, and it is \emph{this} category which is shown to be equivalent to a nice subclass of $\cartSh(X_\bullet^\cover)$ as mentioned in the paragraph above.
This category, denoted by $\Coh_\cover(X)$, consists of complexes of sheaves of modules such that their restriction to any open subset $U$ in some (fixed) cover $\cover$ is quasi-isomorphic to a complex of coherent sheaves, i.e.  complexes $\cal{F}^\bullet$ such that $\cal{F}^\bullet|U\in\Coh(U)$ for any $U\in\cover$.
This category is natural in the choice of cover, i.e.  given any refinement $\anothercover\supset\cover$ we obtain a functor $\Coh_\cover(X)\to\Coh_\anothercover(X)$, and so we can consider the homotopy colimit over all covers.
Using the language of \emph{homotopical categories}, we can localise $\Coh_\cover(X)$ to obtain an $(\infty,1)$-category, and similarly for $\CCoh(X)$, and \cite[Lemmas~4.6.5, 4.6.6]{Hos2020a} shows that the natural embedding $\Coh_\cover(X)\hookrightarrow\CCoh(X)$ induces a fully faithful and essentially surjective functor between the homotopy colimit of the localisations of the $\Coh_\cover(X)$ and the localisation of $\CCoh(X)$.
This does not suffice to say anything about whether or not $\Coh(X)$ and $\CCoh(X)$ are actually equivalent in the analytic case, but it shows that the language of sheaves on the Čech nerve might prove useful in the future.

In fact, we can try to further enrich (a nice full subcategory of) $\Sh(X_\bullet^\cover)$ by endowing its objects with a simplicial version of \emph{Koszul connections}, giving us access to the language of \emph{Chern--Weil theory} and the ability to talk about characteristic classes (i.e.  \emph{Chern classes}).
This was achieved in, again, \cite{Gre1980}, and translated into the language of $\Sh(X_\bullet^\cover)$ and $(\infty,1)$-categories in, again, \cite{Hos2020a}, but we refer the interested reader to these sources instead of overly repeating ourselves here.
The full story makes heavy use of the twisting cochains mentioned at the end of the previous section.

\subsection{String diagrams of composite endomorphisms}
\label{subsection:certain-string-diagrams}

The aim of this section is to show how (pre)sheaves on simplicial spaces might be used to describe things that consist of ``stuff, relations between subsets of the stuff, and linear orders on those relations''.
Hopefully experts in the relevant fields of application can extract the core ideas and refine the technical details to fit (and this is our excuse for a lower level of rigour in what follows).

A slightly more precise (but still intentionally vague and slightly incorrect) description of the main idea is the following:
\begin{quotation}
  \itshape
  Diagrams as in \cref{figure:example-string-diagram} --- consisting of some fixed objects (represented by wires) and some endomorphisms of tensor products (represented by boxes), which describe ``endomorphisms generated by endomorphisms'' --- are in bijection with lax sections of a certain presheaf on a simplicial space that depends only on the objects and the ``unordered shape'' of the diagram, cf. \cref{figure:section-of-sheaf-on-type skeleton}.
  Furthermore, both the presheaf and the simplicial space are $1$-truncated, i.e.  have only degree~$0$ and degree~$1$ parts.
\end{quotation}

\begin{figure}[h!]
    \centering
    \begin{tikzpicture}
      \pgfmathsetmacro{\padding}{0.25}
      \foreach \x in {0,1,2,3}
        \draw (\x,3) to (\x,0);
      \node [thick,fill=white] at (0,3) {$A$};
      \node [thick,fill=white] at (1,3) {$B$};
      \node [thick,fill=white] at (2,3) {$C$};
      \node [thick,fill=white] at (3,3) {$D$};
      \draw [thick,fill=white] (0-\padding,1.75) rectangle ++(1+2*\padding,0.5);
      \draw [thick,fill=white] (1-\padding,0.75) rectangle ++(2+2*\padding,0.5);
      \node at (0.5,2) {$f$};
      \node at (2,1) {$g$};
    \end{tikzpicture}
    \caption{This diagram describes four objects, $A$, $B$, $C$, and $D$, along with two endomorphisms $f\in\End(A\otimes B)$ and $g\in\End(B\otimes C\otimes D)$, with $f$ being applied \emph{before} $g$, i.e.  the entire diagram describes the composite endomorphism $(\id_A\otimes g)\circ(f\otimes\id_C\otimes\id_D)$.}
    \label{figure:example-string-diagram}
\end{figure}
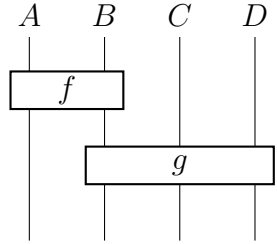

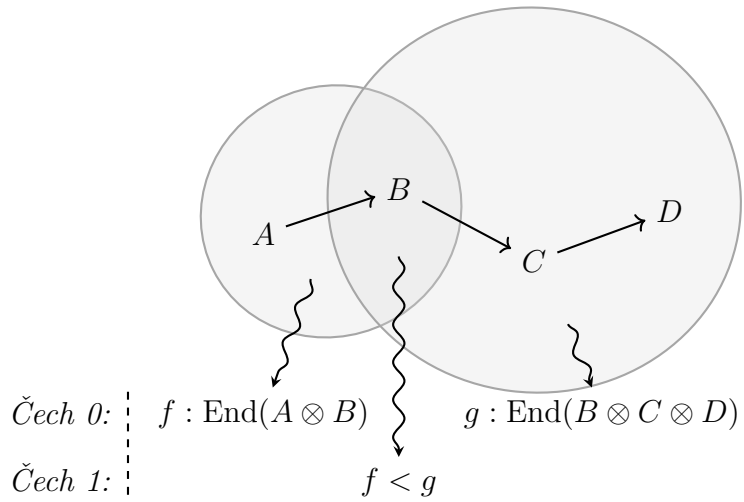
\begin{figure}[h!]
    \centering
    \begin{tikzpicture}[xscale=1.8,yscale=1.2,decoration={snake,segment length=5mm,amplitude=0.75mm}]
      \node [white] (T1) at (0,0) {$A$};
      \node [white] (T2) at (1,0.5) {$B$};
      \node [white] (T3) at (2,-0.3) {$C$};
      \node [white] (T4) at (3,0.25) {$D$};
      \ellipsebyfoci{thick,draw=gray!70,fill=gray!30,fill opacity=0.25}{T1}{T2}{3.5cm};
      \ellipsebyfoci{thick,draw=gray!70,fill=gray!30,fill opacity=0.25}{T2}{T4}{5.5cm};
      % redraw nodes over the ellipses
      \node at (T1) {$A$};
      \node at (T2) {$B$};
      \node at (T3) {$C$};
      \node at (T4) {$D$};
      \draw [thick,->] (T1) to (T2);
      \draw [thick,->] (T2) to (T3);
      \draw [thick,->] (T3) to (T4);
      \draw[thick,dashed] (-1,-1.75) to (-1,-3);
      \node at (-1.5,-1.95) {\emph{Čech 0:}};
      \node (E1e) at (0,-2) {$f:\End(A\otimes B)$};
      \node (E2e) at (2.5,-2) {$g:\End(B\otimes C\otimes D)$};
      \node at (-1.5,-2.7) {\emph{Čech 1:}};
      \node (E1E2o) at (1,-2.75) {$f<g$};
      \draw [thick,-stealth,decorate] (0.35,-0.5) to (E1e);
      \draw [thick,-stealth,decorate] (2.25,-1) to (E2e);
      \draw [thick,-stealth,decorate] (1,-0.25) to (E1E2o);
    \end{tikzpicture}
    \caption{The data of a diagram as in \cref{figure:example-string-diagram} is described by picking an endomorphism on each $0$-intersection in the type skeleton, and an order of these endomorphism on each $1$-intersection in the type skeleton.}
    \label{figure:section-of-sheaf-on-type skeleton}
\end{figure}

Before constructing the presheaf in question and formalising this main idea, we start by breaking down the diagram in \cref{figure:example-string-diagram} into more disparate components, explaining how we arrive at \cref{figure:section-of-sheaf-on-type skeleton}.

Our very first step is to separate out the \emph{shape} of the diagram from the \emph{contents} of the diagram, as shown in \cref{figure:example-string-diagram-skeleton}.
With this new representation, we can now ask ourselves what role \emph{order} is playing.
That is, in the diagram on the left in \cref{figure:example-string-diagram-skeleton},
\begin{enumerate}[(a)]
  \item the order of the wires (i.e.  objects) is important; and
  \item the order of the boxes (i.e.  endomorphisms) is important.
\end{enumerate}
But note that the first one is somehow more of a simple bookkeeping tool, whereas the second actually encodes meaning --- (horizontally) reordering the wires doesn't change the diagram, as long as we reorder the boxes along with them; (vertically) reordering the boxes \emph{does} change the resulting composite endomorphism, in a way that can't be compensated for by relabelling the objects.\footnote{Of course, one can argue that it can be compensated for by simply relabelling the \emph{morphisms}, and ask why we are treating objects as somehow more important than morphisms? This is a very reasonable critique, and emphasises the ``problem'' that we are really trying to solve here: ``\emph{given a tensor product of objects}, how can we describe endomorphisms generated by other endomorphisms?''}
Because of this, it seems like separating out these two uses of order might also be advantageous.
In doing so, we can actually also do away with the string diagram notation altogether: we end up with two groups of data, namely the \emph{type skeleton}, and the \emph{processes}, as shown in \cref{figure:example-string-diagram-fully-separated}.

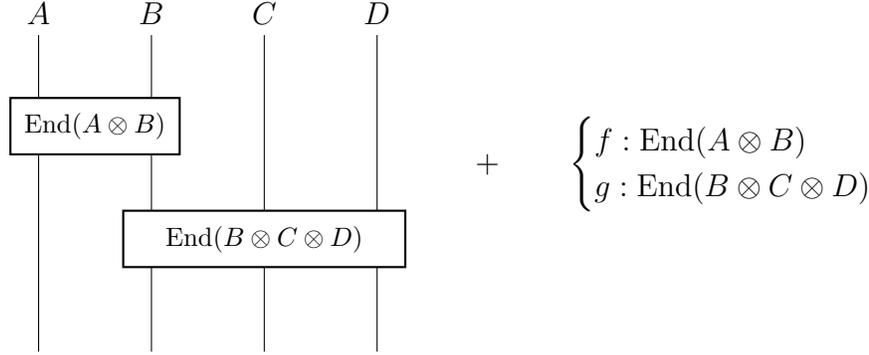
\begin{figure}[h!]
    \centering
    \[
      \begin{tikzpicture}[scale=1.5,baseline={(current bounding box.center)}]
        \pgfmathsetmacro{\padding}{0.25}
        \foreach \x in {0,1,2,3}
          \draw (\x,3) to (\x,0);
        \node [thick,fill=white] at (0,3) {$A$};
        \node [thick,fill=white] at (1,3) {$B$};
        \node [thick,fill=white] at (2,3) {$C$};
        \node [thick,fill=white] at (3,3) {$D$};
        \draw [thick,fill=white] (0-\padding,1.75) rectangle ++(1+2*\padding,0.5);
        \draw [thick,fill=white] (1-\padding,0.75) rectangle ++(2+2*\padding,0.5);
        \node at (0.5,2) {\footnotesize$\End(A\otimes B)$};
        \node at (2,1) {\footnotesize$\End(B\otimes C\otimes D)$};
      \end{tikzpicture}
      \qquad+\qquad
      \begin{cases}
        f:\End(A\otimes B)
      \\g:\End(B\otimes C\otimes D)
      \end{cases}
    \]
    \caption{The diagram from \cref{figure:example-string-diagram}, but where we \emph{start to} separate out the \emph{contents} of the diagram (i.e.  the specific endomorphisms) from the \emph{shape} (i.e.  objects $A$, $B$, $C$, and $D$, and endomorphisms on two specific tensor products of these).}
    \label{figure:example-string-diagram-skeleton}
\end{figure}

\begin{figure}[h!]
    \centering
    \[
      \underbrace{%
        \begin{cases}
          A < B < C < D
        \\A\otimes B,\quad B\otimes C\otimes D
        \end{cases}
      }_{\text{type skeleton}}
      \quad+\quad
      \underbrace{%
        \begin{cases}
          f:\End(A\otimes B)
        \\g:\End(B\otimes C\otimes D)
        \end{cases}
        \quad+\quad
        f < g
      }_{\text{processes}}
    \]
    \caption{The diagram from \cref{figure:example-string-diagram,figure:example-string-diagram-skeleton}, but where we separate out the data into two groups: the \emph{type skeleton} tells us that we have four objects, tells us how to refer to them (i.e.  gives them a linear order), and tells us which ones are ``related'' via a box/endomorphism; the \emph{processes} are the endomorphisms which are applied to our objects, along with the information of the order in which these applications happen.}
    \label{figure:example-string-diagram-fully-separated}
\end{figure}

Our next step is to turn the type skeleton into something that looks more like a topological space, but still describes the same data.
The linear order $A<B<C<D$ gives us a poset, and the data of the types of endomorphisms gives us what we can think of as a basis for a topology on the totally ordered set, as shown in \cref{figure:type-skeleton-as-space}.
Indeed, this actually defines a ($1$-truncated) semi-simplicial\footnote{For simplicity, we finally use \emph{semi}-simplicial instead of simplicial constructions here: the degenerate intersections $U_{\alpha\alpha}$ in the Čech nerve offer no extra information in this setting.} space analogous to the Čech nerve: we obtain $S_\bullet$ given by
\begin{itemize}
  \item $S_0=U\sqcup V :=\{A,B\}\sqcup\{B,C,D\}$;
  \item $S_1=U\cap V=\{B\}$;
  \item $S_p=\varnothing$ for $p\geq2$;
  \item $S_\bullet f_1^i\colon S_1\to S_0$ (for $i=0,1$) is just the inclusion $U\cap V$ into either $U$ or $V$.
\end{itemize}

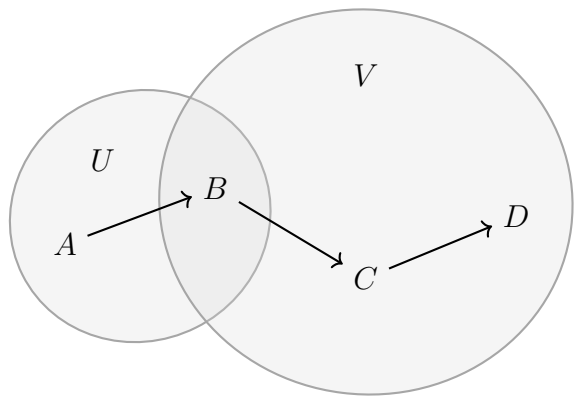
\begin{figure}[h!]
    \centering
    \begin{tikzpicture}[xscale=2,yscale=1.5]
      \node [white] (T1) at (0,0) {$A$};
      \node [white] (T2) at (1,0.5) {$B$};
      \node [white] (T3) at (2,-0.3) {$C$};
      \node [white] (T4) at (3,0.25) {$D$};
      \ellipsebyfoci{thick,draw=gray!70,fill=gray!30,fill opacity=0.25}{T1}{T2}{3.5cm};
      \ellipsebyfoci{thick,draw=gray!70,fill=gray!30,fill opacity=0.25}{T2}{T4}{5.5cm};
      % redraw nodes over the ellipses
      \node at (T1) {$A$};
      \node at (T2) {$B$};
      \node at (T3) {$C$};
      \node at (T4) {$D$};
      \draw [thick,->] (T1) to (T2);
      \draw [thick,->] (T2) to (T3);
      \draw [thick,->] (T3) to (T4);
      % label sets
      \node at (0.25,0.75) {$U$};
      \node at (2,1.5) {$V$};
    \end{tikzpicture}
    \caption{The type skeleton from \cref{figure:example-string-diagram-fully-separated} as a totally ordered set with a basis for a topology (or, indeed, as a semi-simplicial space $S_\bullet$). The two ``open'' sets contain the information of the types of endomorphisms: the fact that $U$ contains $A$ and $B$ relays the information that we have one endomorphism on $A\otimes B$, and similarly for $V$ and $B\otimes C\otimes D$.}
    \label{figure:type-skeleton-as-space}
\end{figure}

Now we can start to build a presheaf $\scr{L}^\bullet\in\Sh(S_\bullet)$ on this type skeleton.
We define
\[
  \scr{L}^0(x_1<\ldots<x_m)
   := \End(x_1<\ldots<x_m)
\]
which, in our specific setting, is just
\[
  \begin{aligned}
    \scr{L}^0(U)
    &= \End(A\otimes B)
  \\\scr{L}^0(V)
    &= \End(B\otimes C\otimes D)
  \end{aligned}
\]
and $\scr{L}^0(S_0)=\scr{L}^0(U)\sqcup\scr{L}^0(V)$.
Next we define
\begin{equation*}
\resizebox{.98\hsize}{!}{$
     \scr{L}^1(x_1<\ldots<x_m)
   := F\big(
  \big\{
    f\in\End(y_1<\ldots<y_n) \mid \{x_1<\ldots<x_m\}\hookrightarrow \{y_1<\ldots<y_n\}
  \big\}
  \big)$},
\end{equation*}
where $F$ takes a set $X$ and returns the $1$-truncated simplicial space with the elements of $X$ as $0$-simplices, and a unique $1$-simplex connecting every (ordered) pair of $0$-simplices;
in our setting, this gives
\begin{equation*}
\resizebox{.98\hsize}{!}{$
  \scr{L}^1(U\cap V)
  =
  \begin{cases}
    f\in\bigsqcup_{k=1}^4\left\{
      \End(x_{i_1}\otimes\ldots\otimes x_{i_k})
      \mid
      x_{i_j}=B\mbox{ for some }1\leq j\leq k
    \right\}
    &\mbox{as $0$-simplices,}
  \\f\to g\quad\mbox{for all $f,g$}
    &\mbox{as $1$-simplices,}
  \end{cases}
$}
\end{equation*}
where the face maps from the $1$-simplices to the $0$-simplices are given by the domain and codomain maps, i.e.  $(f\to g)\mapsto f$ and $(f\to g)\mapsto g$.

Finally then, we can ask what exactly a lax (global) section of $\scr{L}^\bullet$ is.
By \cref{subsection:weak-sections}, using the definition of a section as a map from the (constant sheaf of the) ($1$-truncated) tower of simplices $\Delta^0\hookrightarrow\Delta^1$, this should be
\begin{itemize}
  \item a $0$-simplex $s^0\in\scr{L}^0(S_0)$; and
  \item a $1$-simplex $s^1\in\scr{L}^1(S_1)$
\end{itemize}
such that
\begin{itemize}
  \item $f_1^0\colon s^1\mapsto s^0|U$; and
  \item $f_1^1\colon s^1\mapsto s^0|V$.
\end{itemize}
But, by construction, this is exactly the following:
\begin{itemize}
  \item $s^0_U\in\End(A\otimes B)$ and $s^0_V\in\End(B\otimes C\otimes D)$; and
  \item $s^1\in\{s^0_U\to s^0_V\}\sqcup\{s^0_V\to s^0_U\}$
\end{itemize}
i.e.  an endomorphism $f\in\End(A\otimes B)$, an endomorphism $g\in\End(B\otimes C\otimes D)$, and a choice of order $f<g$ or $g<f$.
That is, exactly the data of a \emph{relation} (in the sense of \cref{figure:example-string-diagram-fully-separated}).

\begin{remark}
  Although we have constructed here a presheaf $\scr{L}^\bullet$ with values in simplicial sets (or topological spaces, graphs, categories, posets, \ldots, depending on how you want to view the ``thing'' that has endomorphisms as ``points'' and a unique directed ``path'' between each ordered pair of points), and then considered lax sections, we could have instead opted to define $\scr{L}^1$ as simply being the set of all formal symbols $f<g$ for all ordered pairs $(f,g)$, and then considered a section in the non-lax sense, i.e.  as a map from $\constant{\{*\}}$ instead of from $\constant{\Delta^0}\hookrightarrow\constant{\Delta^1}$.
\end{remark}

We do not want the underlying (relatively simple) ideas to be lost in the notation and technical construction: all we have really done is constructed an object that encodes the information of ``a choice of element on each open subset, and a choice of order of these elements on each intersection'', as shown in \cref{figure:section-of-sheaf-on-type skeleton}.

\begin{remark}
\label{remark:why-a-sheaf}
  The question as to whether or not the presheaf thus constructed is actually a sheaf can be understood as a question concerning the endomorphisms available to us, i.e.  the objects comprising its (degree-wise) sections.
  That is, the sheaf condition, when unravelled, is equivalent to the following:
  \begin{quotation}
    Any endomorphism $f\in\End(x_1\otimes\ldots\otimes x_n)$ is determined uniquely and entirely by the data of all of its ``partial evaluations'', given by tracing over a subset $\{x_{i_1},\ldots,x_{i_k}\}$ of its arguments.
  \end{quotation}
  That is, we can ensure that the sheaf condition is automatically satisfied for \emph{any} presheaf by imposing some condition on \emph{all} the endomorphisms themselves, i.e.  on the underlying (symmetric) monoidal category.
  Whether or not this is a reasonable thing to do is not a question that we intend to answer here.

  As a side note, since we are \emph{not} endowing our type skeleton with the \emph{Alexandrov} (or \emph{upper-set}) topology, we cannot appeal to the fact that any presheaf on a poset with the Alexandrov topology is automatically a sheaf.
\end{remark}

\begin{remark}
  In the classical world of sheaves, an important construction is the \emph{espace étalé}, which realises a sheaf $\cal{F}\in\Sh(X)$ as a bundle living over $X$.
  We can try to construct something analogous here: consider the simplicial space $L_\bullet$ whose $0$-simplices are all the possible endomorphisms (of all possible types, i.e.  of $T_i$ for all $i$, and of $T_i\otimes T_j$ for all $i,j$, and \ldots), and which has a $1$-simplex between each (ordered) pair of objects;
  then a section of $\cal{L}^\bullet$ is the same as a morphism of simplicial spaces $\sigma\colon S_\bullet\to L_\bullet$ \emph{over $|S_\bullet|$}, i.e.  in our example, such that $\sigma(U)$ consists of endomorphisms $f\in\End(A\otimes B)$ and $\sigma(V)$ consists of endomorphisms of $g\in\End(B\otimes C\otimes D)$.

  We leave the details of this construction vague and incomplete, but we mention it anyway in the hope that some interested reader might try to formalise this definition.
\end{remark}

%%% BIBLIOGRAPHY

\EditInfo{November 23, 2022}{January 17, 2023}{Pasha Zusmanovich}

\end{paper}